\documentclass[11pt]{amsart}
\usepackage[all]{xy}
\usepackage{amsmath}
\usepackage{amsfonts}
\usepackage{amssymb}
\usepackage{latexsym}
\usepackage{graphicx}
\usepackage{color}
\usepackage{amscd}
\usepackage{epsfig}
\usepackage{cite}

\newtheorem{problem}{Problem}[section]
\newtheorem{definition}[problem]{Definition}
\newtheorem{lemma}[problem]{Lemma}
\newtheorem{theorem}[problem]{Theorem}

\title{A sup-Hodge bound for exponential sums}
\author{Chunlei Liu}
\address{Department of Mathematics, Shanghai Jiao Tong
University, Shanghai 200240, P.R. China.} \email{clliu@sjtu.edu.cn}
\thanks{This research is supported by NSFC Grant No.
10671015.}
\begin{document}
\maketitle
\begin{abstract}
The $C$-function of $T$-adic exponential sums is studeid. An
explicit arithmetic bound is established for the Newton polygon of
the $C$-function. This polygon lies above the Hodge polygon. It
gives a sup-Hodge bound of the $C$-function of $p$-power order
exponential sums.
\end{abstract}



\section{Introduction}
Let $p$ be a prime number, $\mathbb{F}_p=\mathbb{Z}/(p)$,
$\overline{\mathbb{F}}_p$ a fixed algebraic closure of
$\mathbb{F}_p$, and $\mathbb{F}_{p^k}$ the subfield of
$\overline{\mathbb{F}}_p$ with $p^k$ elements.

Let $q>1$ be a power of $p$, $W$ the ring scheme of Witt vectors,
$\mathbb{Z}_{q}=W(\mathbb{F}_q)$, $\mathbb{Q}_q$ the fraction field
of $\mathbb{Z}_q$,
$\overline{\mathbb{Q}}_p=\lim\limits_{\stackrel{\rightarrow}{k}}\mathbb{Q}_{p^k}$,
and $\mathbb{C}_p$ the $p$-adic completion of
$\overline{\mathbb{Q}}_p$.

Let $\triangle\supsetneq\{0\}$ be an integral convex polytope
 in $\mathbb{R}^n$, and $I$ the set of vertices of
 $\triangle$ different from the origin.
Let
$$f(x)=\sum\limits_{u\in \triangle}(a_ux^u,0,0,\cdots)\in
W(\mathbb{F}_q[x_1^{\pm1},\cdots,x_n^{\pm1}])\text{ with }
\prod_{u\in I}a_u\neq0,$$ where $x^u=x_1^{u_1}\cdots x_n^{u_n}$ if
$u=(u_1,\cdots,u_n)\in\mathbb{Z}^n$.\begin{definition} For
$k\in\mathbb{N}$, the sum
$$S_{f}(k,T)=\sum\limits_{x\in(\mathbb{F}_{q^k}^{\times})^n}
(1+T)^{{\rm
Tr}_{\mathbb{Z}_{q^k}/\mathbb{Z}_p}(f(x))}\in\mathbb{Z}_p[[T]]$$ is
call a $T$-adic exponential sum. And the function
$$L_f(s,T)=\exp(\sum\limits_{k=1}^{\infty}S_f(k,T)\frac{s^k}{k})\in
1+s\mathbb{Z}_p[[T]][[s]],$$ as a power series in the single
variable $s$ with coefficients in the $T$-adic complete field
$\mathbb{Q}_p((T))$, is called an $L$-function of $T$-adic
exponential sums. \end{definition}

Let $m$ be a positive integer, $\zeta_{p^m}$ a primitive $p^m$-th
root of unity, and $\pi_m=\zeta_{p^m}-1$.  Then $S_f(k,\pi_m)$ is
the exponential sum studied by Liu-Wei \cite{LWe}. If $m=1$, the
exponential sum $S_f(k,\pi_m)$ was studied by
Adolphson-Sperber\cite{AS, AS2}. And, if $n=1$, the exponential sum
$S_f(k,\pi_m)$ was studied by Kumar-Helleseth-Calderbank\cite{KHC}
and Li \cite{Li}.

\begin{definition}The function
$$C_f(s,T) =C_f(s,T;\mathbb{F}_q) =\exp(\sum\limits_{k=1}^{\infty}-(q^k-1)^{-n}S_{f}(k,T)\frac{s^k}{k}),$$
 as a power series in the single
variable $s$ with coefficients in the $T$-adic complete field
$\mathbb{Q}_p((T))$, is called a $C$-function of $T$-adic
exponential sums.\end{definition} We have
$$L_f(s,T) = \prod_{i=0}^n C_f(q^is, T)^{(-1)^{n-i+1}{n\choose i}},$$
and $$C_f(s, T)= \prod_{j=0}^{\infty} L_f(q^js,
T)^{(-1)^{n-1}{n+j-1\choose j}}.$$ So the $C$-function $C_f(s,T)$
and the $L$-function $L_f(s,T)$ determine each other. From the last
identity, one sees that
$$C_f(s,T)\in
1+s\mathbb{Z}_p[[T]][[s]].$$

Let $C(\triangle)$ be the cone generated by $\triangle$, and
$M(\triangle)=M(\triangle)\cap \mathbb{Z}^n$. There is a degree
function $\deg$ on $C(\triangle)$ which is
$\mathbb{R}_{\geq0}$-linear and takes the values $1$ on each
co-dimension $1$ face not containing $0$. For $a\not\in
C(\triangle)$, we define $\deg(a)=+\infty$.
\begin{definition}A convex function on
$[0,+\infty]$ which is linear between consecutive integers with
initial value $0$ is called the infinite Hodge polygon of
$\triangle$ if its slopes between consecutive integers are the
numbers $\deg(a)$, $a\in M(\triangle)$. We denote this polygon by
$H_{\triangle}^{\infty}$.
\end{definition}

Liu-Wan \cite{LWa} also proved the following.

\begin{theorem}[Hodge bound]We have $$T-\text{adic NP of }C_{f}(s,T)\geq\text{ord}_p(q)(p-1)H_{\triangle}^{\infty},$$
where NP is the short for Newton polygon.\end{theorem} Denote by
$\lceil x\rceil$ the least integer equal or greater than $x$, and by
$\{x\}$ the fractional part of $x$.
\begin{definition}Let $C\subseteq M(\triangle)$ be a finite subset. We define
$$r_C=\max_{\beta}(\#\{a\in C\mid
\{\deg(pa)\}'\geq\beta\}-\#\{a\in C\mid \{\deg(a)\}'\geq\beta\}).$$
\end{definition}
\begin{definition}Let $a\subseteq M(\triangle)$. We define
$$\varpi(a)=\lceil p\deg(a)\rceil-\lceil\deg(a)\rceil+r_{\{a\in M(\triangle)\mid
\deg(u)<\deg(a)\}\cup\{a\}}-r_{\{a\in M(\triangle)\mid
\deg(u)<\deg(a)\}}.$$
\end{definition}
\begin{definition}The arithmetic polygon $p_{\triangle}$
of $\triangle$ is a convex function on $[0,+\infty]$ which is linear
between consecutive integers with initial value $0$, and whose
slopes between consecutive integers are the numbers
$\varpi_{\triangle}(a)$, $a\in M(\triangle)$.
\end{definition}

One can prove the following.
\begin{theorem}\label{arithhodge}We have
then$$p_{\triangle}\geq (p-1)H_{\triangle}^{\infty}.$$ Moreover,
they coincide at the point $n!{\rm Vol}(\triangle)$.
\end{theorem}

Let $D$ be the least positive integer such that
$\deg(M(\triangle))\subseteq\frac{1}{D}\mathbb{Z}$. The main result
of this paper is the following.
\begin{theorem}\label{main1}If $p>3D$, then
$$T-\text{adic NP of }C_{f}(s,T)\geq {\rm ord}_p(q)p_{\triangle}.$$
\end{theorem}
 From the above theorem we shall deduce
the following.
\begin{theorem}\label{main2}If $p>3D$, then, for $t\in\mathbb{C}_p$ with $0\neq|t|_p<1$, we have
$$t-\text{adic NP of }C_{f}(s,t)\geq{\rm ord}_p(q)p_{\triangle}.$$
\end{theorem}
\section{The $T$-adic Dwork Theory}
In this section we review the $T$-adic analogue of Dwork theory on
exponential sums.

Let
$$E(t)=\exp(\sum_{i=0}^{\infty}\frac{t^{p^i}}{p^i})=\sum\limits_{i=0}^{+\infty}\lambda_it^i
\in 1+t{\mathbb Z}_p[[t]]$$ be the $p$-adic Artin-Hasse exponential
series. Define a new $T$-adic uniformizer $\pi$ of ${\mathbb
Q}_p((T))$ by the formula $E(\pi)=1+T$. Let $\pi^{1/D}$ be a fixed
$D$-th root of $\pi$. Let
$$L=\{\sum_{u\in M(\triangle)}c_u\pi^{\deg(u)}x^u:\
 c_u\in\mathbb{Z}_q[[\pi^{1/D}]] \}.$$

Let $a\mapsto\hat{a}$  be the Teichm\"{u}ller lifting. One can show
that the series
$$E_f(x) :=\prod\limits_{a_u\neq0}E(\pi \hat{a}_ux^u)\in L.$$
Note that the Galois group of $\mathbb{Q}_q$ over $\mathbb{Q}_p$ can
act on $L$ but keeping $\pi^{1/D}$ as well as the variable $x$
fixed.
 Let $\sigma$ be the Frobenius element in the
Galois group such that $\sigma(\zeta)=\zeta^p$ if $\zeta$ is a
$(q-1)$-th root of unity. Let $\Psi_p$ be the operator on $L$
defined by the formula
$$\Psi_p(
 \sum\limits_{i\in
M(\triangle)} c_ix^i)=\sum\limits_{i\in M(\triangle)} c_{pi}x^i.$$
Then $\Psi:=\sigma^{-1}\circ\Psi_p\circ E_f$ acts on the $T$-adic
Banach module
$$B=\{\sum\limits_{u\in M(\triangle)}c_i\pi^{\deg(u)}x^u \in L,\
 \text{\rm ord}_T(c_u)\rightarrow+\infty
 \text{ if }\deg(u)\rightarrow+\infty\}.$$
We call it Dwork's $T$-adic semi-linear operator because it is
semi-linear over $\mathbb{Z}_q[[\pi^{\frac{1}{D}}]]$ .

Let $b=\log_pq$. Then the $b$-iterate $\Psi^b$ is linear over
$\mathbb{Z}_q[[\pi^{1/D}]]$, since
$$\Psi^{b}=\Psi_p^{b}\circ
\prod\limits_{i=0}^{b-1}E_{f}^{\sigma^i}(x^{p^i}).$$ One can show
that $\Psi$ is completely continuous in the sense of Serre
\cite{Se}. So $\det(1-\Psi^bs\mid
B/\mathbb{Z}_{q}[[\pi^{\frac{1}{D}}]])$ and $\det(1-\Psi s\mid
B/\mathbb{Z}_p[[\pi^{\frac{1}{D}}]])$ are well-defined.

We now state the $T$-adic Dwork trace formula\cite{LWa}.
\begin{theorem}
We have
$$C_f(s,T)=\det(1-\Psi^bs\mid
B/\mathbb{Z}_{q}[[\pi^{\frac{1}{D}}]]).$$
\end{theorem}
\begin{lemma}\label{q2p}The Newton polygon of
 $\det(1-\Psi^bs^b\mid
B/\mathbb{Z}_{q}[[\pi^{\frac{1}{D}}]])$ coincides with that of
$\det(1-\Psi s\mid B/\mathbb{Z}_p[[\pi^{\frac{1}{D}}]])$.\end{lemma}
\proof Note that
$$\det(1-\Psi^b s\mid B/\mathbb{Z}_p[[\pi^{\frac{1}{D}}]])
={\rm Norm}(\det(1-\Psi^b s\mid
B/\mathbb{Z}_q[[\pi^{\frac{1}{D}}]])),$$ where Norm is the norm map
from $\mathbb{Z}_q[[\pi^{\frac{1}{D}}]]$ to
$\mathbb{Z}_p[[\pi^{\frac{1}{D}}]]$. The lemma now follows from the
equality
$$\prod\limits_{\zeta^b=1}\det(1-\Psi\zeta s\mid B/\mathbb{Z}_p[[\pi^{\frac{1}{D}}]])
=\det(1-\Psi^b s^b\mid B/\mathbb{Z}_p[[\pi^{\frac{1}{D}}]]).$$
\endproof

Write
$$\det(1-\Psi s\mid B/\mathbb{Z}_p[[\pi^{\frac{1}{D}}]])=\sum\limits_{i=0}^{+\infty}(-1)^ic_is^i.$$
\begin{theorem}\label{coefficient-q2p}
The $T$-adic Newton polygon of
 $\det(1-\Psi^bs\mid
B/\mathbb{Z}_{q}[[\pi^{\frac{1}{D}}]])$ is the lower convex closure
of the points
$$(m,\text{ord}_{T}(c_{bm})),\ m=0,1,\cdots.$$\end{theorem}

\proof By Lemma \ref{q2p}, the $T$-adic Newton polygon of
 $\det(1-\Psi^bs^b\mid
B/\mathbb{Z}_{q}[[\pi^{\frac{1}{D}}]])$ is the lower convex closure
of the points
$$(i,\text{ord}_{T}(c_i)),\ i=0,1,\cdots.$$
It is clear that $(i,\text{ord}_{T}(c_i))$ is not a vertex of that
polygon if $b\nmid i$. So that Newton polygon is the lower convex
closure of the points
$$(bm,\text{ord}_{T}(c_{bm})),\ m=0,1,\cdots.$$
It follows that the $T$-adic Newton polygon of
 $\det(1-\Psi^bs\mid
B/\mathbb{Z}_{q}[[\pi^{\frac{1}{D}}]])$ is the lower convex closure
of the points
$$(m,\text{ord}_{T}(c_{bm})),\ m=0,1,\cdots.$$
\endproof
\section{The arithmetic bound}
In this section we prove the following.
\begin{theorem}We have
$${\rm ord}_T(c_{bm})\geq p_{\triangle}(m).$$
\end{theorem}
\proof First, we choose a basis of
$B\otimes_{\mathbb{Z}_p}\mathbb{Q}_p(\pi^{1/D})$ over
$\mathbb{Q}_p(\pi^{1/D})$ as follows. Fix a normal basis
$\bar{\xi}_i$, $i\in\mathbb{Z}/(b)$ of $\mathbb{F}_q$ over
$\mathbb{F}_p$. Let $\xi_i$ be their Terchm\"{u}ller lift of
$\bar{\xi}_i$. The system $\xi_i$, $i\in\mathbb{Z}/(b)$ is a normal
basis of $\mathbb{Q}_q$ over $\mathbb{Q}_p$. Then
$\{\xi_ix^u\}_{u\in M(\triangle),1\leq i\leq b}$ is a basis of
$B\otimes_{\mathbb{Z}_p}\mathbb{Q}_p(\pi^{1/D})$ over
$\mathbb{Q}_p(\pi^{1/D})$.

Secondly, we write out the matrix of $\Psi$ on
$B\otimes_{\mathbb{Z}_p}\mathbb{Q}_p(\pi^{1/D})$ with respect to the
basis $\{\xi_ix^u\}_{u\in M(\triangle),1\leq i\leq b}$. Write
$$E_f(x)= \sum\limits_{u\in M(\triangle)}\gamma_ux^u,$$
and
$$\sigma^{-1}(\xi_j\gamma_{pu-w})=\sum\limits_{i=1}^{b}\gamma_{(u,i),(w,j)}\xi_i.$$
Then
$$\Psi(\xi_jx^w)=\sum\limits_{u\in
M(\triangle)}\sigma^{-1}(\xi_j\gamma_u)\Psi_p(x^{u+w})$$$$=\sum\limits_{u\in
M(\triangle)}\sigma^{-1}(\xi_j\gamma_{pu-w})x^u$$$$=\sum\limits_{u\in
M(\triangle)}\sum\limits_{i=1}^{b}\gamma_{(u,i),(w,j)}\xi_ix^u.
$$
So $(\gamma_{(u,i),(w,j)})_{u,w\in M(\triangle),1\leq i,j\leq b}$ is
the matrix of $\Psi$ on
$B\otimes_{\mathbb{Z}_p}\mathbb{Q}_p(\pi^{1/D})$ with respect to the
basis $\{\xi_ix^u\}_{u\in M(\triangle),1\leq i\leq b}$.

Thirdly, we claim that$${\rm
ord}_T(\gamma_{(u,i),(w,j)})\geq\lceil\deg(pu-w)\rceil.$$ In fact,
this follows from the
equality$$\sigma^{-1}(\xi_j\gamma_{pu-w})=\sum\limits_{i=1}^{b}\gamma_{(u,i),(w,j)}\xi_i,$$
and the inequality ${\rm ord}_T(\gamma_u)\geq\lceil\deg(u)\rceil$.

Finally, we show that
$${\rm ord}_T(c_{bm})\geq p_{\triangle}(m).$$
Note
that
$$c_{bm}=\sum\limits_{A}\text{det}((\gamma_{(u,i),(w,j)})_{(u,i),(w,j)\in
A}),
$$
where $A$ runs over all subsets of $
M(\triangle)\times\mathbb{Z}/(b)$ with cardinality $bm$. So it
suffices to show that
$${\rm ord}_T(\det(\gamma_{(i,u),(j,\omega)})_{(i,u),(j,\omega)\in
A})\geq bp_{\triangle}(m).$$ Note that
$$\det(\gamma_{(i,u),(j,\omega)})_{(i,u),(j,\omega)\in A})
=\sum\limits_{\tau\in S_A}\sum\limits_{a\in A}{\rm
ord}_{\pi}(\gamma_{a,\tau(a)}),$$ where $S_A$ is the permutation
group of $A$. So it suffices to show that
$$\sum\limits_{a\in A}{\rm
ord}_{\pi}(\gamma_{a,\tau(a)})\geq bp_{\triangle}(m),\ \tau\in
S_A.$$ Since $${\rm
ord}_T(\gamma_{(u,i),(w,j)})\geq\lceil\deg(pu-w)\rceil,$$ the
theorem follows from the following. \endproof
\begin{theorem}If $p>3D$, $A$ is a subset of $
M(\triangle)\times\mathbb{Z}/(b)$ with cardinality $bm$, and
$\tau\in S_A$, then$$\sum\limits_{a\in
A}\lceil\deg(p\nu(a)-\nu(\tau(a)))\rceil\geq bp_{\triangle}(m),$$
where $\nu(u,i)=u$.\end{theorem}\proof
 We have
$$\sum\limits_{a\in
A}\lceil\deg(p\nu(a)-\nu(\tau(a)))\rceil\geq \sum\limits_{a\in
A}\lceil\deg(p\nu(a))-\nu(\tau(a))\rceil$$$$\geq \sum\limits_{a\in
A}\lceil\deg(p\nu(a))\rceil-\lceil\deg(\nu(\tau(a)))\rceil+1_{\{\deg(p\nu(a))\}'>\{\deg(\nu(\tau(a)))\}'}$$
$$\geq\sum\limits_{a\in
A}\lceil\deg(p\nu(a))\rceil-\lceil\deg(\nu(a))\rceil+1_{\{\deg(p\nu(a))\}'>\{\deg(\nu(\tau(a)))\}'}.$$
 Choose a set $B$ of cardinality $|A|$ such that $B\cap A$ is as big as possible under the condition
that, for some $\alpha$,
$$\{a\in M(\triangle)\times\mathbb{Z}/(b)\mid \deg(\nu(a))< \alpha\}\subseteq B\subsetneq\{a\in
M(\triangle)\times\mathbb{Z}/(b)\mid \deg(\nu(a))\leq \alpha\}.
$$
Choose a permutation $\tau_0$ on $B\cap A$ which agrees with $\tau$
on $(B\cap A)\cap\tau^{-1}(B\cap A)$. Extend it trivially to $B$. We
have
$$\sum\limits_{a\in
A}\lceil\deg(p\nu(a))\rceil-\lceil\deg(\nu(a))\rceil
\geq\sum\limits_{a\in
B}(\lceil\deg(p\nu(a))\rceil-\lceil\deg(\nu(a))\rceil)+2\#(A\setminus
B).$$ We also have$$\sum\limits_{a\in
A}1_{\{\deg(p\nu(a))\}'>\{\deg(\nu(\tau(a)))\}'}\geq\sum\limits_{a\in
B}1_{\{\deg(p\nu(a))\}'>\{\deg(\nu(\tau_0(a)))\}'}-2\#(A\setminus
B).$$ It follows that
$$\sum\limits_{a\in
A}\lceil\deg(p\nu(a))\rceil-\lceil\deg(\nu(a))\rceil+1_{\{\deg(p\nu(a))\}'>\{\deg(\nu(\tau(a)))\}'}$$
$$\geq\sum\limits_{a\in
B}\lceil\deg(p\nu(a))\rceil-\lceil\deg(\nu(a))\rceil+1_{\{\deg(p\nu(a))\}'>\{\deg(\nu(\tau_0(a)))\}'}$$
$$\geq\sum\limits_{a\in
B}(\lceil\deg(p\nu(a))\rceil-\lceil\deg(\nu(a))\rceil)+r_B,$$ where
$$r_B=\max_{\beta}(\#\{a\in B\mid
\{\deg(p\nu(a))\}'\geq\beta\}-\#\{a\in B\mid
\{\deg(\nu(a))\}'\geq\beta\}).$$ Choose a set $C$ of cardinality $m$
such that for some $\alpha$,
$$\{a\in M(\triangle)\mid \deg(\nu(a))< \alpha\}\subseteq C\subsetneq\{a\in
M(\triangle)\mid \deg(\nu(a))\leq \alpha\}.$$ Recall that
$$r_C=\max_{\beta}(\#\{a\in C\mid \{\deg(pa)\}'\geq\beta\}-\#\{a\in
C\mid \{\deg(a)\}'\geq\beta\}).$$ It is easy to see that $r_B=br_C$,
and $$\sum\limits_{a\in
B}(\lceil\deg(p\nu(a))\rceil-\lceil\deg(\nu(a))\rceil)+r_B$$$$=b\sum\limits_{a\in
C}(\lceil\deg(pa)\rceil-\lceil\deg(a)\rceil)+br_C=bp_{\triangle}(m).$$
The theorem now follows.
\endproof

\end{document}